\documentclass[graybox]{svmult}

\usepackage{mathptmx}       
\usepackage{helvet}         
\usepackage{courier}        
\usepackage{type1cm}        
\usepackage{makeidx}         
\usepackage{graphicx}        
\usepackage{multicol}        
\usepackage[bottom]{footmisc}

\usepackage{amsmath,amssymb,amscd}
\usepackage{amsfonts}
\usepackage{grffile}
\usepackage[nooneline,hang]{subfigure}
\usepackage[normalem]{ulem}
\usepackage{rotating}

\newcommand{\ie}{i.e., }
\newcommand{\eg}{e.g., }
\newcommand{\lift}{\mathcal{L}}
\newcommand{\tskip}{t_\mathrm{skip}}

\newcommand{\restrict}{\mathcal{R}}

\newcommand{\avg}[1]{\langle #1 \rangle}

\newcommand{\R}{\mathbb{R}}

\renewcommand{\phi}{\varphi}

\renewcommand{\epsilon}{\varepsilon}

\newcommand{\Fig}[1]{Fig.~\ref{#1}}

\definecolor{darkgreen}{rgb}{0.0,0.6,0.0}

\makeindex             

\begin{document}

\title*{Equation-Free Analysis of Macroscopic Behavior in Traffic and Pedestrian Flow}
\author{Christian Marschler, Jan Sieber, Poul G. Hjorth and Jens Starke}
\institute{Christian Marschler (\email{chrms@dtu.dk}), Poul G. Hjorth (\email{pghj@dtu.dk}),
  Jens Starke (\email{jsta@dtu.dk}) \at Department of Applied Mathematics and Computer Science, Technical
University of Denmark, DK-2800 Kongens Lyngby,  Denmark, %
\and Jan Sieber (\email{j.sieber@exeter.ac.uk})\at College of Engineering, Mathematics and Physical
Sciences, University of Exeter, EX4 4QF Exeter, United Kingdom}%
%
%
\maketitle

\abstract{Equation-free methods 
make possible an analysis of the evolution 
of a few coarse-grained or macroscopic quantities for a detailed and
realistic model with a large
number of fine-grained or microscopic variables, 
  even though no
  equations are explicitly given on the macroscopic level. This will
  facilitate a study of how the model behavior depends on parameter
  values including an understanding of transitions between different
  types of qualitative behavior. These methods are introduced and
  explained for traffic jam formation and emergence of oscillatory
  pedestrian counter flow in a corridor with a narrow door.}


\section{Introduction}
\label{sec:intro}
The study of pedestrian and traffic dynamics leads naturally to a
description by a few macroscopic, \eg averaged, quantities of the
systems at hand. On the other hand, so-called microscopic models, \eg
multiagent systems, inherit individual properties of the agents and
can therefore be made very realistic. Among more successful
microscopic models are social force models for pedestrian
dynamics~\cite{Helbing1995,Helbing2001,Corradi2012} and optimal
velocity models in traffic
dynamics~\cite{Bando1995,Gasser2004,Sugiyama2008,Orosz2005,Marschler2013}.
Although computer
simulations of microscopic models for specific scenarios are
straightforward to perform it is often more relevant and useful to
look at the systems on a coarse scale, \eg to investigate a few
macroscopic quantities like first-order moments of distributions or
other macroscopic descriptions which are motivated by the application.

The analysis of the macroscopic behavior of microscopically defined
models is possible by the so-called \emph{equation-free} or
\emph{coarse analysis}.  This 
approach is motivated and
justified by the observation, that multi-scale systems, \eg
many-particle systems, often exhibit low-dimensional behavior. This
concept is well known in physics as slaving of many degrees of freedom
by a few slow variables, sometimes refered to as ``order parameters''
(see e.g.\ \cite{Haken1983a,Haken1983}) and is formalized
mathematically for slow-fast systems by Fenichel's
theory~\cite{Fenichel1979}.  These methods aim for a description of
the system in terms of a small number of variables, which describe the
interesting dynamics. This results in a dimension reduction from many
degrees of freedom to a few degrees of freedom. For example, in
pedestrian flows, we reduce the full system of equations of motion
with equations of motion for each single pedestrian to a
low-dimensional system for weighted mean position and velocity of the
crowd.

A difficulty for such a macroscopic analysis is that governing
equations for the coarse variables, \ie the order parameters, are
often not known. Those equations are often very hard or sometimes even
impossible to derive from first principles especially in models with a
very complicated microscopic dynamics. To extract information about
the macroscopic behavior from the microscopic models equation-free
methods
\cite{kevrekidisgear2003,Kevrekidis2004,KevrekidisSamaey2009,Kevrekidis2010}
can be used. This is done by using a special scheme for switching
between microscopic and macroscopic levels by restriction and lifting
operators and suitably initialized short microscopic simulation bursts
in between. Problems with the initialization of the microscopic
dynamics, \ie the so-called lifting error, have been studied
in~\cite{Marschler2013}. An \emph{implicit equation-free method} for
simplifying the lifting procedure has been introduced, allowing for
avoiding lifting errors up to an error which can be estimated for
reliable results~\cite{Marschler2013}. The equation-free methodology
is most suitable in cases where governing equations for coarse
variables are either not known, or when one wants to study finite-size
effects if the number of particles is too large for investigation of
the full system, but not large enough for a
continuum limit. 
It is even possible to
apply equation-free and related techniques in experiments, where the
microscopic simulation is replaced by observations of an
experiment~\cite{sieberkrauskopf08,Bureau2013,Barton2013}.

For pedestrian and for traffic problems, a particularly interesting case
is a systematic study of the influence of parameters on solutions of
the system. This leads to equation-free bifurcation analysis. One
obtains qualitative as well as quantitative information about the
solutions and their stability. Furthermore, it saves computational time and is
therefore advantageous over a brute-force analysis or computation. The
knowledge of parameter dependence and the basin of attraction of
solutions is crucial for controling systems and ensuring their
robustness. Changes of solutions are summarized in bifurcation diagrams and solution
branches are usually obtained by means of numerical
continuation. These techniques from numerical bifurcation analysis can be combined with
equation-free methods to gain insight into the macroscopic behavior in
a semi-automatic fashion.

In the following, we apply \emph{equation-free bifurcation analysis}
to two selected problems in traffic and pedestrian dynamics. Section
\ref{sec:eqfree} gives a short overview about equation-free
methods. The methods introduced in Section~\ref{sec:eqfree} are
then applied to study traffic jams in the optimal velocity model
(cf.~\cite{Bando1995,Marschler2013}) in Section
\ref{sec:traffic}. Section \ref{sec:peds} describes the macroscopic
analysis of two pedestrian groups in counterflow through a bottleneck
(cf.~\cite{Corradi2012}) and Section \ref{sec:conc} concludes the
paper with a brief discussion and an outlook on future research
directions.


\section{Equation-Free Methods}
\label{sec:eqfree}
\begin{figure}[t]
  \centering
\includegraphics{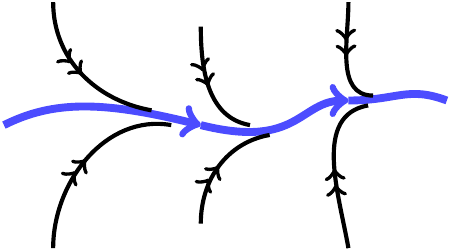}
\caption{Fast convergence to a slow manifold (thick blue
  curve). Trajectories in many dynamical systems converge very
  quickly to a slow manifold, along which the long-time macroscopic behavior takes
  place.}
  \label{fig:slowfast}
\end{figure}

Equation-free methods have been introduced
(cf.~\cite{KevrekidisSamaey2009,Kevrekidis2010} for reviews) to study
the dynamics of multi-scale systems on a macroscopic level without the
need for an explicit derivation of macroscopic equations from the
microscopic model. The necessary information is obtained by
suitably initialized short simulation bursts of the microscopic
system at hand. Equation-free methods assume that the system under
investigation can be usefully described on a coarse scale. Evolution
equations on the macroscopic level are not given explicitly. A big
class of suitable systems are slow-fast systems, which have a
separation of time scales. Under quite general assumptions
(cf.~\cite{Fenichel1979}) these systems quickly converge to a
low-dimensional object in phase space, the so-called \emph{slow
  manifold} (cf. Fig. \ref{fig:slowfast}). The long-term dynamics
(\ie the macroscopic behavior) happens on this slow manifold,
which is usually of much lower dimension than the overall phase space
(of the microscopic system). The goal of equation-free methods is to
gain insight into the dynamics on this slow manifold.

In the following we discuss the equation-free methodology in
detail. The construction of a so-called macroscopic time stepper
requires three ingredients to be provided by the user: the lifting
$\lift$ and restriction $\restrict$ operators to communicate between
the microscopic and macroscopic levels and vice versa, and the
microscopic time stepper $M$. Due to a separation of time scales, it
is possible to construct the macroscopic time stepper by a
\emph{lift-evolve-restrict-scheme}.  This scheme is subsequently used
to perform bifurcation analysis and numerical continuation.

\subsection*{Microscopic time stepper $M$}
To be specific, let us consider a microscopic model in the form of
  a high-dimensional system of $N$ differential equations
\begin{equation}
  \label{eq:micro}
  \dot{u} = f(u)\mbox{.}
\end{equation}
This can be any model of traffic or pedestrian dynamics, possibly
depending on a set of parameters. We generally assume that the number
of degrees of freedom and thereby the dimension $N$ of $u$ is large.
Note that a second-order model, \eg the social force model with forces
$f_\mathrm{force}(x)$, can be written as a first-order model of the
type \eqref{eq:micro} by including the velocities $\dot x = v$ into
the equation. Then $u$ has the form $u=(x,v)$, and the right-hand side
is $f(u)=f((x,v))=[v,f_\mathrm{force}(x)]$. We assume that a
microscopic time stepper $M$ for model \eqref{eq:micro} is
available. That is, we have a routine $M$ (usually a simulation or
software package) with two inputs: the time $t\in\R$ by which we want
to evolve and the initial state $u_0\in\R^N$ from which we start. The output
$M(t,u_0)\in\R^N$ is defined by the relation
\begin{equation}
  \label{eq:microts}
  u(t_0+t) = M(t,u(t_0))\mbox{.}
\end{equation}
That is $M(t,u_0)$ is the state $u$ of \eqref{eq:micro} after
time $t$, starting from $u_0$ at time $t_0$.

\subsection*{Separation of time scales}
We also assume that the dynamics on the macroscopic scale can be
described by a few macroscopic variables $x\in\R^n$, where $n$ is much
smaller than the phase space dimension $N$ of the microscopic model.
This assumption is typically true in many-particle systems, \eg
pedestrian flow and traffic problems. The goal of equation-free
methods is then to construct a time stepper for $x$ on the macroscopic level,
\begin{equation}
  \label{eq:macrots}
  x(t_0+t) = \Phi(t,x(t_0))\mbox{,}
\end{equation}
based on repeated and appropriately initialized runs, \ie simulation
bursts, of the microscopic time stepper $M$ for $u$. In practice, a
user of equation-free methods begins with the identification of a map,
the so-called restriction operator
\begin{displaymath}
  \restrict:\R^N\to\R^n\mbox{,}
\end{displaymath}
which reduces a given microscopic state $u\in\R^N$ to a value of the desired
macroscopic variable $x\in\R^n$. The assumption about the variables $x$
describing the dynamics at the macroscopic scale has to be made more
precise. We require that for all relevant initial conditions $u$ and a
sufficiently long transient time $\tskip$ the result of the microscopic time
stepper \eqref{eq:microts} is (at least locally and up to a small error) uniquely determined by its
restriction, \ie its macroscopic behavior. That is, if for two initial conditions $u_0$ and $u_1$
the relation
\begin{equation}
  \label{eq:rtrans}
  \begin{split}
    &\restrict M(\tskip;u_0)=\restrict M(\tskip;u_1)\quad\mbox{holds, then}\\
    &|\restrict M(\tskip+t;u_0)-\restrict
    M(\tskip+t;u_1)|<C\exp(\epsilon t-\gamma\, \tskip)
  \end{split}
\end{equation}
for all $t\geq0$. In \eqref{eq:rtrans} the pre-factor $C$ should be of
order unity and independent of the choice of $t$, $u_0$ and $u_1$.
The growth rate $\epsilon$ is also assumed to be
smaller than the decay rate $\gamma$. This is what we refer to as
\emph{separation of time scales} between macroscopic and microscopic
dynamics. Requirement \eqref{eq:rtrans} makes the statement ``the
dynamics of $u$ on long time scales can be described by the
macroscopic variable $x=\restrict
u$'' more precise. We also
see that the error in this description can be made as small as desired
by increasing the \emph{healing time} $\tskip$. In fact, requirement
\eqref{eq:rtrans} determines what a good choice of $\tskip$ is for a
given problem.

In order to complete the construction of the macroscopic time stepper
$\Phi$, the user has to provide a lifting operator
\begin{displaymath}
  \lift:\R^n\to \R^N\mbox{,}
\end{displaymath}
which reconstructs a microscopic state $u$ from a given macroscopic
state $x$.  See \cite{Gear2005,Zagaris2009,Zagaris2012} for proposals
how to construct good lifting operators for explicit equation-free
methods (see Eq. ~\eqref{eq:explicitscheme} below). In the case of
implicit equation-free methods the choice of a lifting operator is not
as delicate \cite{Marschler2013}. Also note that the choice of lifting
operator is not unique.

\subsection*{Macroscopic time stepper $\Phi$}
\begin{figure}[t]
  \centering
\includegraphics[width = .9\textwidth]{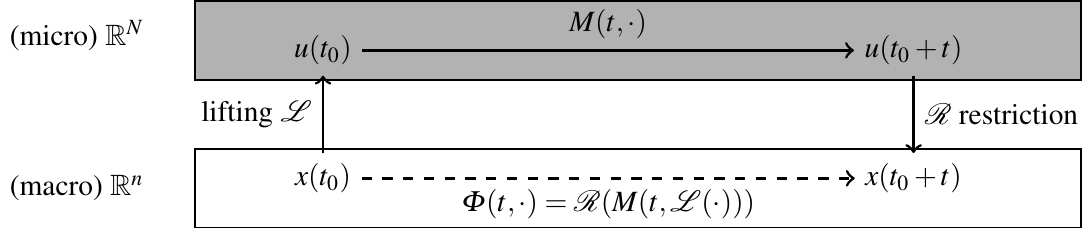} 
  \caption{Scheme for construction of the macroscopic time stepper $\Phi$
    using the lifting $\lift$ and restriction operator $\restrict$ for
  switching between microscopic and macroscopic levels. $M$ denotes
  the microscopic time stepper.}
  \label{fig:eqfree}
\end{figure}

\begin{figure}[t]
  \centering
\includegraphics{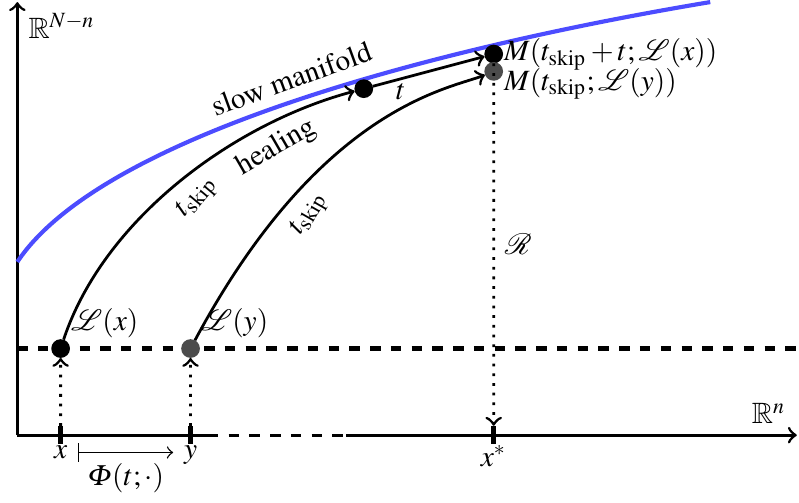}
\caption{Visualization of the implicit scheme
  \eqref{eq:implicitscheme}. The macroscopic time stepper $\Phi$ maps
  the macroscopic state $x$ to the yet unknown macroscopic state
  $y$. The scheme \emph{lift-evolve-restrict} is applied to both
  states. Additionally to the healing step $\tskip$ the dynamics on
  the slow manifold are observed for state $x$ for an additional
  (long) time $t$. Both ``paths'' are compared at the macroscopic end
  point $x^*$. Note, that this scheme defines $y$ implicitly.}
  \label{fig:slowman}
\end{figure}

We can now assemble the approximate macroscopic time stepper $\Phi$
for $x$ by applying the steps \emph{Lift-Evolve-Restrict}, as illustrated in
Fig.~\ref{fig:eqfree} in a judicious manner (cf. Fig.~\ref{fig:slowman}
for a detailed construction): the time-$t$ image
$y=\Phi(t;x)$ of an initial condition $x\in\R^n$ is defined as the
solution $y$ of the implicit equation
\begin{equation}
  \label{eq:implicitscheme}
  \restrict M(\tskip;\lift y)=\restrict M(\tskip+t;\lift x)\mbox{.}
\end{equation}
  Note, that the macroscopic time stepper has originally been introduced
as the explicit definition (cf. also Fig.~\ref{fig:eqfree})
\begin{equation}
  \label{eq:explicitscheme}
  \tilde\Phi(t;x) = \restrict M(t;\lift x).
\end{equation}
The explicit method \eqref{eq:explicitscheme} requires that the
lifting operator maps onto (or very close to) the slow manifold for every macroscopic
point $x$. The implicit method \eqref{eq:implicitscheme} does not have
this requirement and should be the method of choice (cf. the discussion in Section
\ref{sec:traffic}). The implementation of the explicit and implicit
time stepper is further illustrated in Table \ref{tab:algo} using
pseudocode.
Equation \eqref{eq:implicitscheme} is a nonlinear but in
general 
regular system of $n$ equations for the $n$-dimensional variable $y$.
Note that the construction \eqref{eq:implicitscheme} does not require
an explicit derivation of the right-hand side $F:\R^n\to\R^n$ of the
assumed-to-exist macroscopic dynamical system
\begin{equation}
  \label{eq:macro}
  \dot{x} = F(x)\mbox{.}
\end{equation}
However, it can be used to evaluate (approximately) the right-hand
side $F$ in desired arguments $x$ (see below).
The convergence of the time stepper $\Phi$ to the correct time-$t$ map
$\Phi_*$ of the assumed-to-exist macroscopic equation \eqref{eq:macro}
is proven in detail in \cite{Marschler2013}. The error
$|\Phi(t;x)-\Phi_*(t;x)|$ is of order $\exp(\epsilon t-\gamma
\tskip)$.
\begin{table}[t]
  \centering
  \begin{tabular}{|p{\textwidth}|}\hline
    required functions: lift, evolve, restrict (cf. main text)\\
  solution at time t0: x
\begin{verbatim}
function res = Phi(t,x)
u1 = lift(x); u2 = evolve(t,u1); res = restrict(u2);
end
\end{verbatim}\vspace{-.4cm}\\ \hline
  \end{tabular}

  \begin{tabular}{|p{.25\textwidth}|p{.742\textwidth}|}\hline
explicit scheme & implicit scheme\\ \hline
\vspace{-.4cm}
\begin{verbatim}
y = Phi(t,x);
\end{verbatim}
 & \vspace{-.4cm}
\begin{verbatim}
choose dy, tol, y[0] = x, n = 0, err = 2*tol

function res = F(y)
  res = Phi(tskip,y) - Phi(tskip+t,x);
end

while err > tol
  Fy      = F(y[n]); 
  dF      = Jacobian(F,y[n],dy);
  y[n+1]  = y[n] - (dF)^(-1)*(Fy);
  err     = abs(y[n+1] - y[n]); 
  n       = n+1; 
end 
y         = y[n];
\end{verbatim}\vspace{-.5cm}\\
 \hline
  \end{tabular}
  \caption{Pseudocode algorithm for computing the macroscopic solution
    $y$ after time $t$ using the macroscopic time stepper for the
    solution $x$ using the explicit \eqref{eq:explicitscheme} and
    implicit \eqref{eq:implicitscheme} scheme, respectively. The
    implicit scheme uses a Newton iteration with a given tolerance $\texttt{tol}$ to find
    $y$. For one-dimensional $y$ the Jacobian {\tt dF} is given by
    $\texttt{(F(y[n]+dy)-F(y[n]))/dy}$. Note, that the complexity of the
    implicit scheme stems mainly from the Newton iteration, which is
    not specific for equation-free computations.}
  \label{tab:algo}
\end{table}
\subsection*{Advantages of  equation-free methods}
What additional benefits can the macroscopic time stepper $\Phi$ have
beyond simulation of the low-dimensional dynamics (which could have
been accomplished by running long-time simulations using $M$
directly)?
\begin{itemize}
\item \emph{Finding locations of macroscopic equilibria}
  regardless of their dynamical stability: macroscopic equilibria $x$
  are given by solutions to the $n$-dimensional implicit
  equation $\Phi(t_0;x)=x$, or, in terms of lifting and restriction:
  \begin{equation}\label{eq:fixpoint}
    \restrict M(\tskip+t_0;\lift x)=\restrict M(\tskip;\lift x)
  \end{equation}
  for a suitably chosen time $t_0$ (a good choice is of the same order
  of magnitude as $\tskip$). The stability of an equilibrium $x$,
  found by solving \eqref{eq:fixpoint}, is determined by solving the
  generalized eigenvalue problem $Ax=\lambda Bx$ with the Jacobian
  matrices
  \begin{displaymath}
    A=\frac{\partial}{\partial x}\restrict M(\tskip+t_0;\lift x)\mbox{,\quad}
    B=\frac{\partial}{\partial x}\restrict M(\tskip;\lift x)\mbox{.}    
  \end{displaymath}
  Stability is determined by the modulus of the eigenvalues $\lambda$
  (where $|\lambda|<1$ corresponds to stability).
\item \emph{Projective integration} of \eqref{eq:macro}: one can
  integrate the macroscopic system \eqref{eq:macro} by point-wise
  approximation of the right-hand side $F$ and a standard numerical
  integrator. For example, the explicit Euler scheme for
  \eqref{eq:macro} would determine the value $x_{k+1}\approx
  x((k+1)\Delta t)$ from $x_k\approx x(k\Delta t)$ implicitly by approximating
  \begin{displaymath}
    F(x_k)=\frac{1}{\delta}\left[\restrict M(\tskip+\delta;\lift x_k)-
      \restrict M(\tskip;\lift x_k)\right]
  \end{displaymath}
  with a small time $\delta$, and then solving the implicit equation 
  \begin{displaymath}
    \restrict M(\tskip;\lift x_{k+1})-\restrict M(\tskip;\lift x_k)=F(x_k)
  \end{displaymath}
with respect to $x_{k+1}$.
  Projective integration is useful if the macroscopic time step
  $\Delta t$ can be chosen such that $\Delta t\gg \delta$, or for
  negative $\Delta t$, enabling integration backward in time for the
  macroscopic system \eqref{eq:macro}.
\item \emph{Matching the restriction}: Sometimes it is useful to find a
  ``realistic'' microscopic state $u$, corresponding to a given
  macroscopic value $x$. ``Realistic'' corresponds in this context to
  ``after rapid transients have settled''. This can be
  accomplished by solving the nonlinear equation
  \begin{equation}
    \label{eq:match}
    \restrict M(\tskip;\lift y)=x
  \end{equation}
  for $y$ and then setting $u=M(\tskip;\lift y)$.
\end{itemize}
The formulas \eqref{eq:fixpoint} and \eqref{eq:match} have already
been presented and tested in \cite{Vandekerckhove2011}, where they
were found to have vastly superior performance compared to alternative
proposals for consistent lifting (such as presented in
\cite{Gear2005,Zagaris2009,Zagaris2012}).

\subsection*{Bifurcation analysis and numerical continuation}
Building on top of the basic uses of the macroscopic time stepper
$\Phi$, one can also use advanced tools for the study of
parameter-dependent systems. Suppose that the microscopic time stepper
$M$ (and, thus, the macroscopic time stepper $\Phi$) depends on a
system parameter $p$. We are interested in how macroscopic equilibria
and their stability change as we vary $p$. In the examples in Sections
\ref{sec:traffic} and \ref{sec:peds} the primary system
parameter is the target velocity (traffic) and door width
(pedestrians), respectively. 

\begin{figure}[t]
  \centering
  \includegraphics[width = .625\textwidth]{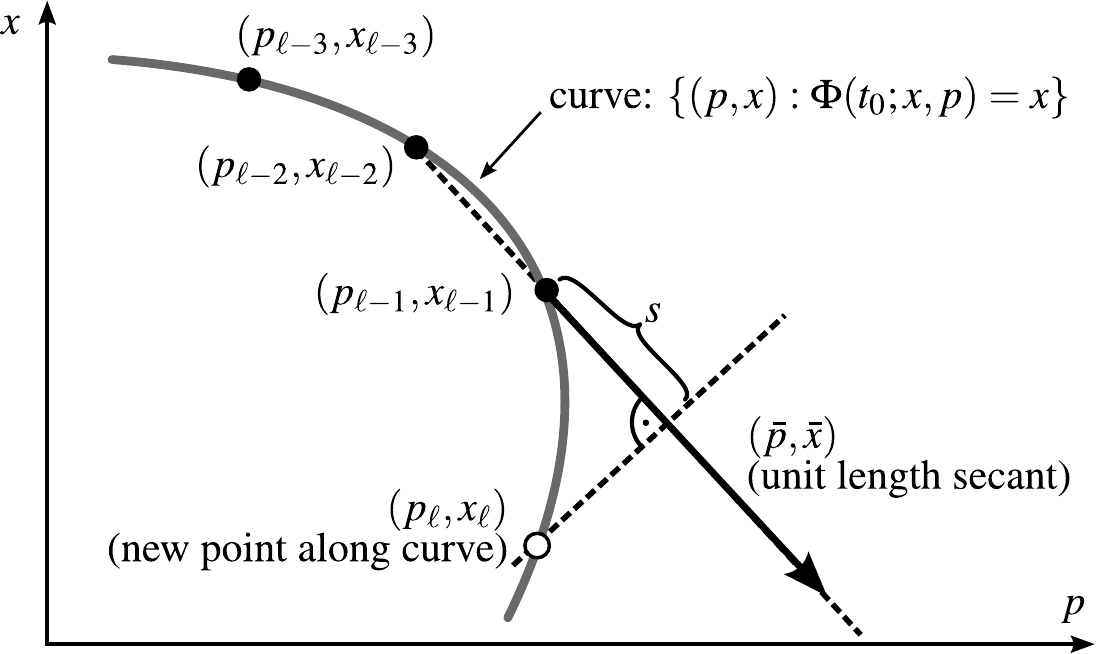}
  \caption{Pseudo-arclength continuation of a curve of fixed points
    $\lbrace (p,x): \Phi(t_0;x,p) = x\rbrace$ of the macroscopic time
    stepper $\Phi$. A new point $(\bar p, \bar x)$ is computed along
    the secant through $(p_{\ell-2},x_{\ell-2})$ and
    $(p_{\ell-1},x_{\ell-1})$ in a so-called predictor step. The
    following corrector step solves the equilibrium condition (cf. \eqref{eq:pcont}) in the
    perpendicular direction to find the next equilibrium $(p_\ell,x_\ell)$ on the curve.}
  \label{fig:pseudoarc}
\end{figure}
When tracking equilibria in a parameter-dependent problem one may
start at a parameter value $p_0$, where the desired equilibrium $x_0$
(given by $\Phi(t_0;x_0,p_0)=x_0$) is stable so that it can be found
by direct simulations. This achieves a good initial guess, which is
required to solve the nonlinear equations \eqref{eq:fixpoint} reliably
with a Newton iteration for near-by $p$ close to $p_0$. In the traffic
system studied in Section \ref{sec:traffic} the equilibrium
corresponding to a single phantom jam undergoes a saddle-node
bifurcation (also called fold, that is, the equilibrium turns back in
the parameter changing its stability, see Fig.~\ref{fig:liftcomp}(a)
for an illustration). In order to track equilibria near folds one
needs to extend the nonlinear equation for the macroscopic equilibrium
with a so-called pseudo-arclength condition, and solve for the
equilibrium $x$ and the parameter $p$ simultaneously
\cite{Beyn2002,Kuznetsov2004}. That is, suppose we have already found
a sequence $(p_k,x_k)$, $k=1,\ldots,\ell-1$, of equilibria and
parameter values. We then determine the next pair $(p_\ell,x_\ell)$ by
solving the extended system for $(p_\ell,x_\ell)$:
\begin{equation}\label{eq:pcont}
  \begin{aligned}
    0&=\Phi(t_0;x_\ell,p_\ell)-x_\ell\quad&\mbox{equilibrium condition}\\
    s&=\bar p_\ell(p_\ell-p_{\ell-1})+
    \bar x^T_\ell(x_\ell-x_{\ell-1})\quad&\mbox{pseudo-arclength condition.}\\
    \end{aligned}
\end{equation}
The vector
\begin{equation}
  \label{eq:secant}
   (\bar p_\ell,\bar
    x^T_\ell)=\frac{(p_{\ell-1}-p_{\ell-2},x^T_{\ell-1}-x^T_{\ell-2})}{
      |(p_{\ell-1}-p_{\ell-2},x^T_{\ell-1}-x^T_{\ell-2})|}
 \end{equation}
 is the secant through the previous two points, scaled to unit length,
 and $s$ is the approximate desired distance of the newly found point
 $(p_\ell,x_\ell)$ from its predecessor $(p_{\ell-1},x_{\ell-1})$. The
 continuation method \eqref{eq:pcont} permits one to track equilibria
 through folds such as shown in Fig.~\ref{fig:liftcomp}(a) or Hopf
 bifurcations such as shown in Fig.~\ref{fig:pedres}(b) (where the
 equilibrium becomes unstable and small-amplitude oscillations
 emerge).  For a more detailed review on methods for bifurcation
 analysis the reader is referred to standard
 references, \eg \cite{Beyn2002,Kuznetsov2004}.
\section{Traffic Models}
\label{sec:traffic}
We apply the methods introduced in Section \ref{sec:eqfree} to the
optimal velocity (OV) model~\cite{Bando1995} as an example of
microscopic traffic models. The model captures the main features
of experiments of cars on a ring road~\cite{Sugiyama2008}. We
exploit equation-free numerical bifurcation analysis to answer the following
questions; 1) for which parameter values in the OV  model do
we expect traffic jams and 2) how severe are they?

The equations of motion for car $n$ in the OV model are
\begin{equation}
  \label{eq:ovode}
  \tau \ddot{x}_n + \dot{x}_n = V(x_{n+1}-x_n), \qquad V(\Delta x_n) = v_0 (\tanh(\Delta x_n-h) + \tanh(h)),
\end{equation}
where $\tau=0.588$ is the reaction time and $V$ is the optimal
velocity function depending on the velocity parameter $v_0$ and
inflection point $h$. Periodic boundary conditions $x_{n+N} = x_n + L$
are used for $N=60$ cars on a ring road of length $L=60$.
Depending on the choice of $v_0$ and $h$ one observes uniform
flow, \ie all cars have \emph{headway} $\Delta x_n = 1$, or a traffic jam,
\ie a region of high density of cars. It is worth noting, that
bistable parameter regimes can exist, \ie a stable uniform flow and a
stable traffic jam coexist and one or the other emerges, depending on initial conditions.

First, we fix $h=1.2$ and study the bifurcation diagram in dependence
of $v_0$. Before we are able to apply the algorithms presented in
Section \ref{sec:eqfree}, we have to define the lifting and
restriction operators.

\subsection*{The restriction and lifting operators}
The restriction operator $\restrict$, used to compute the macroscopic
variable to describe phenomena of interest (here the deviation of the density
profile from a uniform flow) of the microscopic model on a coarse
level, is chosen as the standard deviation of the distribution of
headway values
\begin{equation}
  \label{eq:restraf}
  \restrict (u) = \sigma = \sqrt{\frac{1}{N-1} \sum_{n=1}^N \left(\Delta x_n -\langle \Delta x \rangle\right)^2},
\end{equation}
where $\langle \Delta x \rangle$ is the mean headway.

As the numerical continuation operates in a local neighborhood of the
states, the lifting operator can be based on a previously computed
microscopic reference state $\tilde{u} = (\tilde{x},\tilde{y})$ for
positions $\tilde x$ and velocities $\tilde{y}$ and its macroscopic
image under $\restrict$, $\tilde \sigma = \restrict \tilde u$. We use
$\tilde u$ and $\tilde \sigma$ to obtain a microscopic profile $u$ for
every $\sigma\approx\tilde \sigma$:
\begin{equation}
  \label{eq:lifttraf}
  \lift_\mu (\tilde u, \sigma) = u=(x,y) = \left(x_\text{new}, 
    V(x_\text{new})\right), \quad
  x_\text{new} = \frac{\mu \sigma}{\tilde\sigma}(\Delta \tilde x -\avg{\Delta \tilde x})+\avg{\Delta \tilde x}\mbox{.}
\end{equation}
We let the lifting $\lift_\mu$ depend on an artificial parameter
$\mu$.  We will vary $\mu$ later to demonstrate that the resulting
bifurcation diagram is independent of the particular choice of
$\lift$. 

\subsection*{Numerical results}
\begin{figure}[t]
  \centering \subfigure[{Bifurcation diagram,
    $h=1.2$}]{\label{fig:liftcomp_bif}\includegraphics[width
    =.46\textwidth]{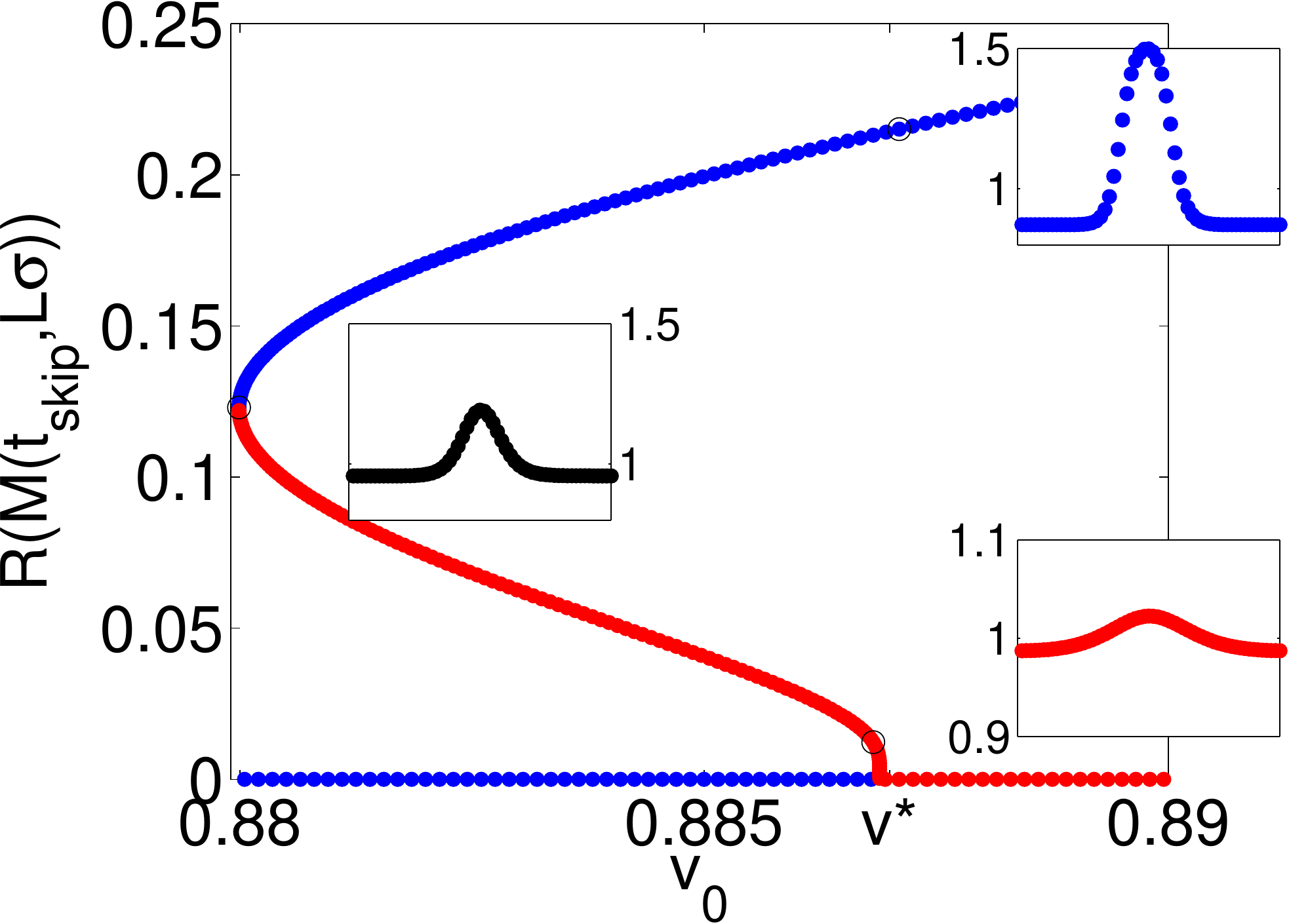}} \hfill
  \subfigure[{Unhealed bifurcation diagrams
    \newline$\mu=0.95,1.00,1.05$}]{\label{fig:liftcomp_dir}\includegraphics[width
    =.45\textwidth]{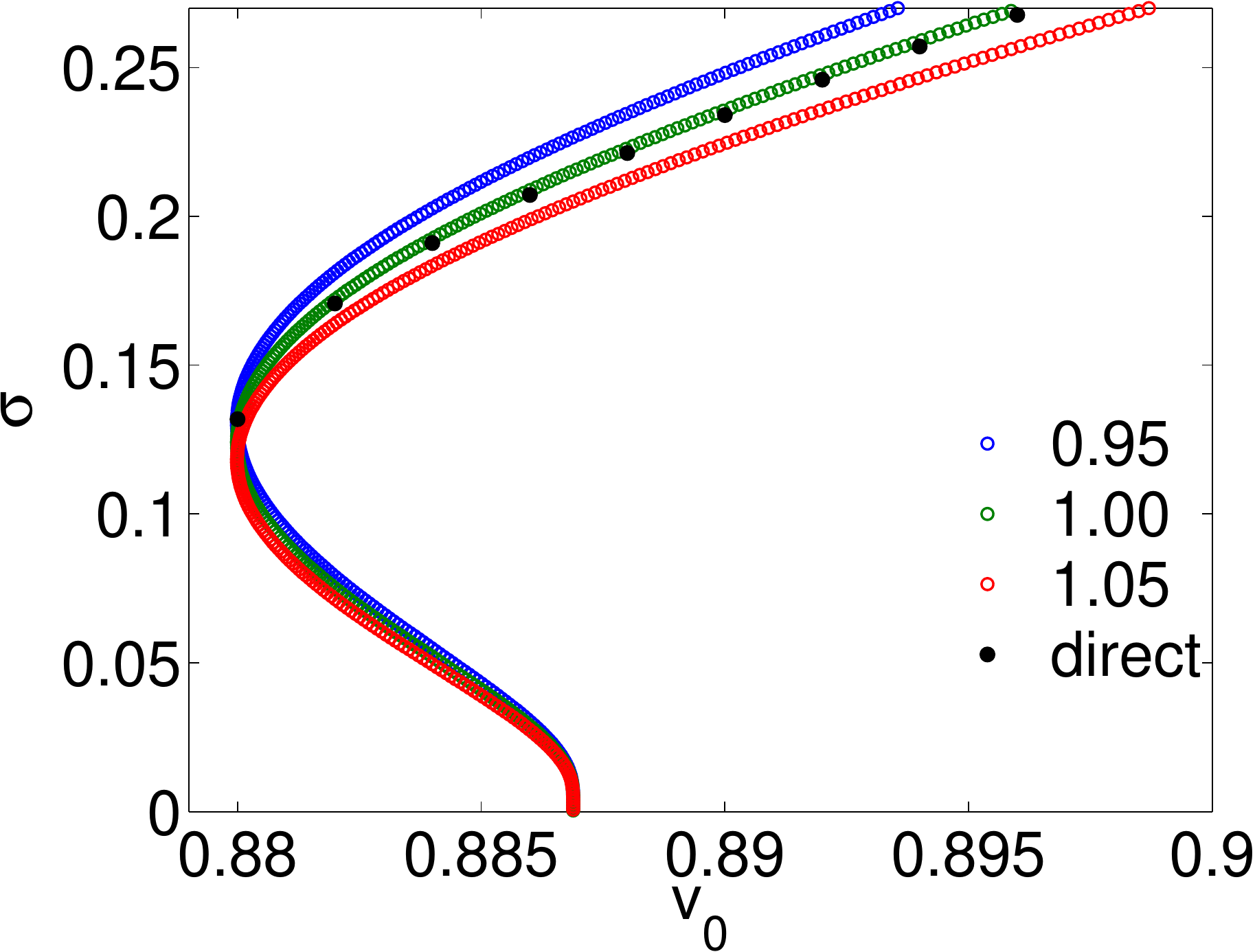}}
  \subfigure[{Two-parameter
    continuation}]{\label{fig:foldcont}%
\setlength{\unitlength}{0.48\textwidth}
    \begin{picture}(1,0.8)
\put(0.4,0.5){\textbf{jam}}
\put(0.6,0.2){\textbf{free flow}}
\put(0.72,0.43){\begin{rotate}{45}\textbf{coexistence}\end{rotate}}
\includegraphics[width=.45\textwidth]{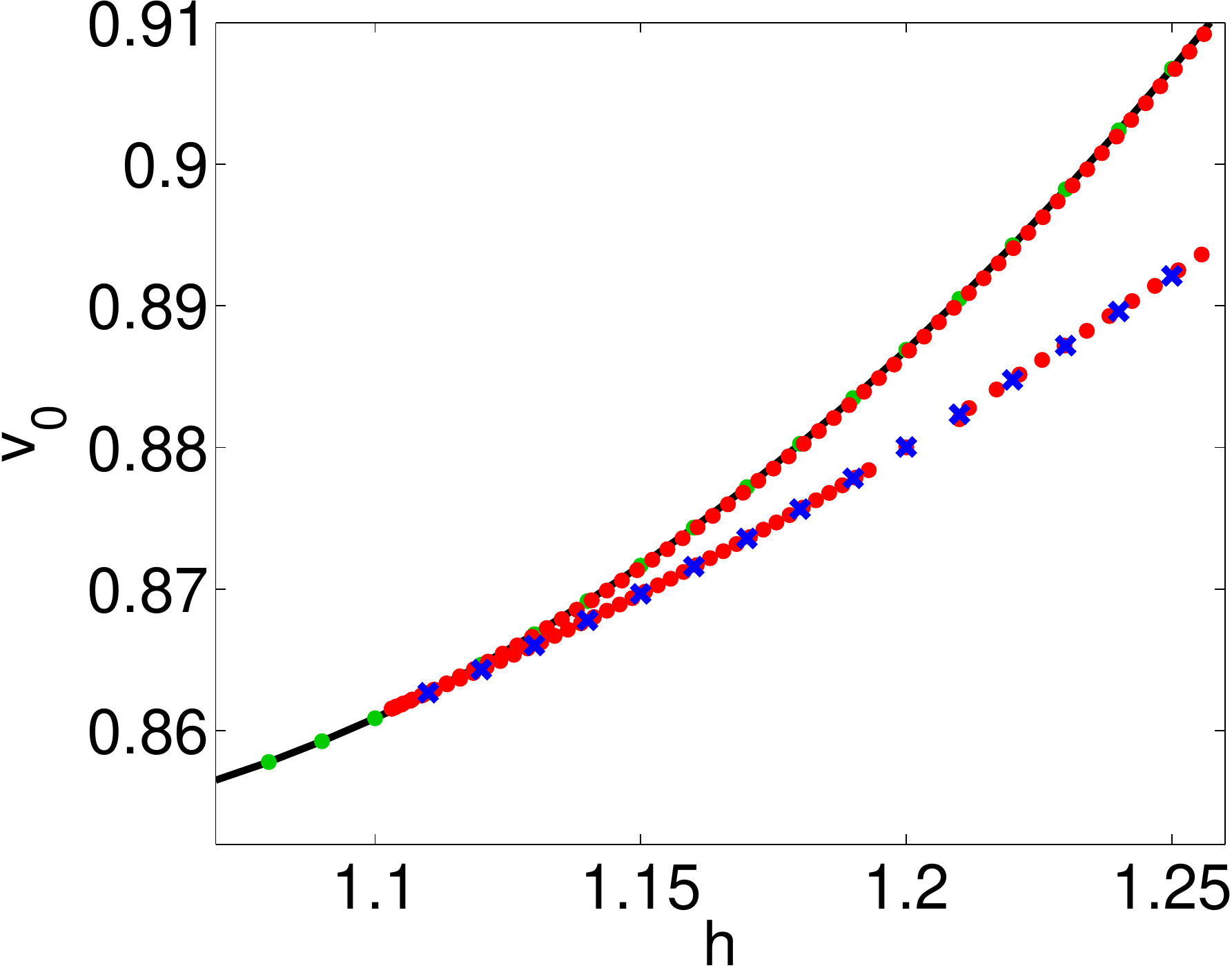}%
    \end{picture}}
 \hfill \subfigure[{Healed
    bifurcation diagrams
    \newline$\mu=0.95,1.00,1.05$}]{\label{fig:liftcomp_dirheal}\includegraphics[width
    =.47\textwidth]{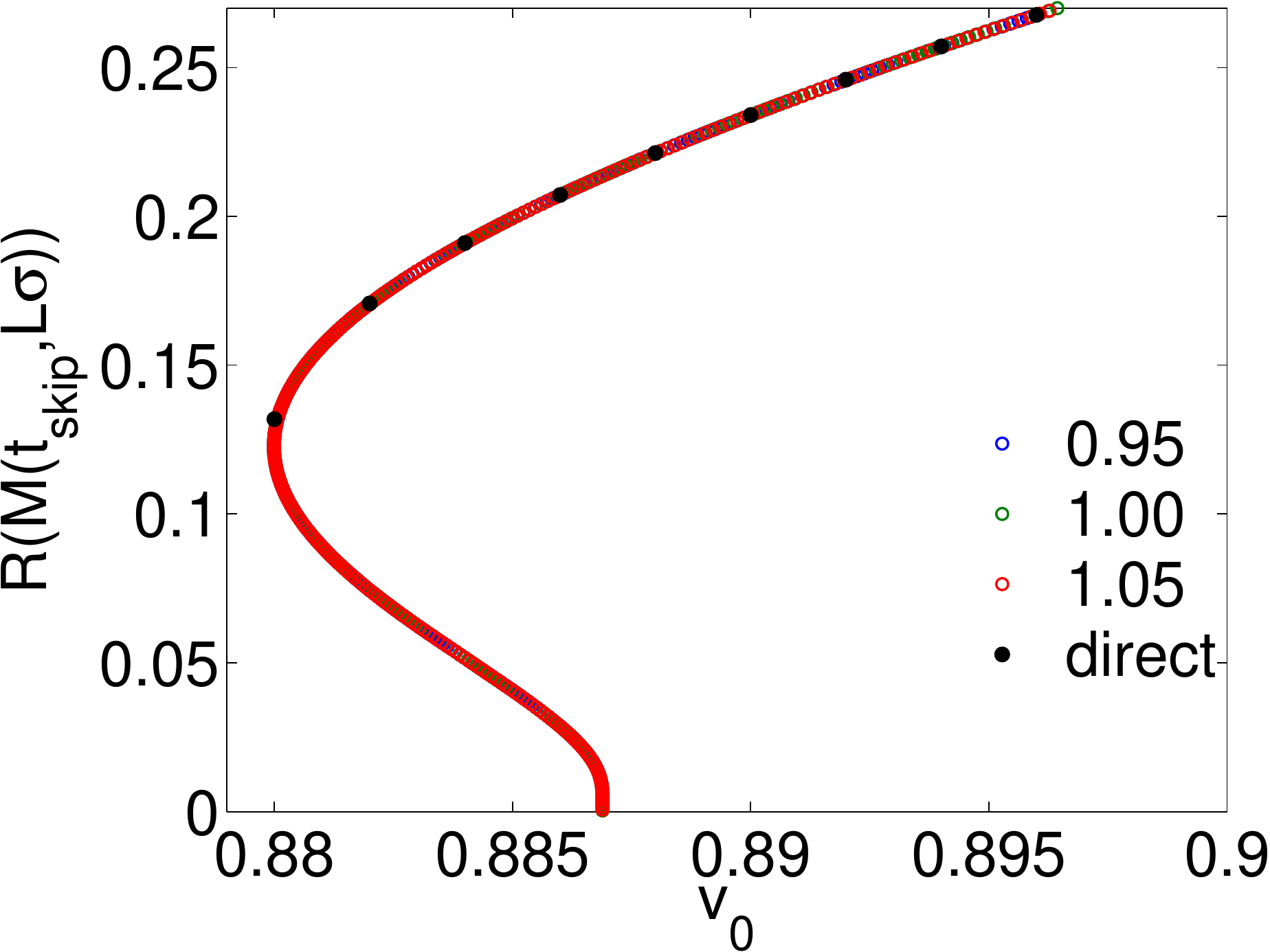}}%
  \caption{Equation-free bifurcation analysis for the optimal velocity
    model \eqref{eq:ovode}. \textup{(a)} Bifurcation diagram in healed quantities for
    $h=1.2$. Headway profiles are shown for selected points (black
    circles) along the branch. Blue and red dots denote stable and
    unstable solutions, respectively. \textup{(b)} and \textup{(d)}
    show bifurcation diagrams for different lifting operators. Healed
    values in \textup{(d)} lie exactly on the same branch and recover
    the results from direct simulation (black dots). Thus, the
      choice of lifting operator $\lift$ does not affect the results
      if one reports the healed values (in contrast to \textup{(b)},
      reporting the solutions $\sigma$ of \eqref{eq:fixpoint}).
    \textup{(c)} Two-parameter bifurcation diagram for continuation of
    the fold point. Saddle-node (blue crosses) and Hopf points (green
    dots) from measurements in one-dimensional diagrams are in perfect
    agreement with the continuation in two parameters $h$ and $v_0$
    and the analytical curve (black line).}
  \label{fig:liftcomp}
\end{figure}
The results of the equation-free bifurcation analysis are shown in
\Fig{fig:liftcomp}. The bifurcation diagram for fixed $h=1.2$ (cf.\
\Fig{fig:liftcomp_bif}) shows a stable traffic jam for parameter
values $v_0 > v^* = 0.887$. By continuation of the solution from a
stable traffic jam towards smaller values of $v_0$ a saddle-node
bifurcation is found at $v_0=0.88$. The traffic jam loses stability
and an unstable solution exists for $v_0 \in [0.88,0.887]$. Continuing
further along the branch, a Hopf bifurcation, \ie a macroscopic
pitchfork bifurcation, where traffic jams are born as small-amplitude
time-periodic patterns, is found at $v_0=0.887$. At this point, stable
uniform flow solutions ($\sigma = 0$) change their stability to
unstable uniform flow solutions. For $v_0 \in [0.88,0.887]$ two stable
solutions coexist. In this one-dimensional system, the unstable
solution separates the stable and the unstable fixed point, acting as
a barrier. Thus, the bifurcation diagram also informs us about the
magnitude of the disturbance
necessary 
  to change the behavior of the system from a stable traffic jam to a
stable free flow. Headway profiles are shown for selected points along
the branch to illustrate the microscopic solutions.
In \Fig{fig:liftcomp_dir} and \Fig{fig:liftcomp_dirheal} the
comparison of different lifting operators is shown. While the
  unhealed values $\sigma$ (cf.\ \Fig{fig:liftcomp_dir}) of the
  equilibrium depends on the choice of $\mu$, the healed values
  $\restrict M(\tskip;\lift\sigma)$, used in the implicit
  equation-free methods (cf.\ \Fig{fig:liftcomp_dirheal} and
  \cite{Marschler2013}) 
  are in perfect agreement with results from direct simulations (black
  dots).

In order to study the dependence on both parameters $v_0$ and $h$
simultaneously, we use an extended set of equations to continue
the saddle-node bifurcation point in \Fig{fig:foldcont}. Blue crosses
and green dots denote measurements of the saddle-node and Hopf
points from one-parameter continuations, respectively. The
two-parameters continuation (red dots) is in perfect agreement with
the measurements. As a check of validity, the Hopf curve (black line
below red dots)
can be computed analytically (cf.~\eg \cite{Marschler2013}) and is shown for comparison.

In conclusion, the analysis pinpoints the parameter values for the
onset and collapse of traffic jams. This information is of potential
use to understand the role of speed limits. The two-parameter
bifurcation diagram in \Fig{fig:foldcont} shows a free flow regime for
small $v_0$ and large $h$ (bottom right part of the diagram). On the
other hand, a large velocity parameter $v_0$ and a small safety distance
$h$ lead to traffic jams (top left part). In between, a coexistence
between free flow and traffic jams is found. The final state depends
on the initial condition. A speed limit lower than the saddle-node
values is necessary to assure a global convergence to the uniform free
flow.

\section{Pedestrian Models}
\label{sec:peds}
For further demonstration of the equation-free bifurcation analysis,
we also apply it to a social force model describing pedestrian
flow~\cite{Helbing1995,Seyfried2006}. A particular setup with two
crowds passing a corridor with bottleneck \cite{Zhang2012} from
opposite sites (the crowd marked blue moving to the right, the crowd
marked red moving to the left) is analyzed with respect to qualitative
changes of the system
behavior~\cite{Corradi2012,MarschlerStarkeLiuKevrekidis14}. To this
end, a coarse bifurcation analysis is used to determine which
bifurcations occur and thereby to understand which solutions are
expected to exist.  Details about the model and the analysis of the
bottleneck problem can be found in~\cite{Corradi2012}. Here, we focus
on the coarse analysis of the problem.

Two parameters have been chosen as the main bifurcations parameters;
the ratio of desired velocities of the two crowds $r_{v_0}=v_0^r/v_0^b$ and
the width of the door $w$ acting as a bottleneck. Microscopic
simulations of the model for two crowds of size $N=100$ reveal two
fundamentally different regimes of the dynamics. One finds a blocked
state and a state that is oscillating at the macroscopic level
(cf.\ \Fig{fig:snaps}) for small and large door widths,
respectively. The question we would like to answer is: how and
  where does the transition from a blocked to an oscillating state
happen? In mathematical terms the question is, where is the
bifurcation point and what type of
bifurcation is observed at the transition? 

\subsection*{The restriction and lifting operators}
We define the macroscopic quantity $m$ as
\begin{equation}
  \label{eq:pedmacro}
  m = \frac{m_r + m_b}{2}, \qquad m_{(r,b)} = \frac{\sum_{i\in(r,b)} \kappa(x_i)
    x_i}{\sum_{i\in(r,b)} \kappa(x_i)},
\end{equation}
where $m_{(r,b)}$ is a weighted average of the longitudinal component
for the blue and red pedestrian crowd, respectively. $\kappa$ gives
more weight to pedestrians close to the door (see
\cite{Corradi2012} for details). Since we expect
oscillations from microscopic observations the pair of variables
$(m,\dot m)$ is used as the macroscopic variable for the
equation-free methods. The transient from the initial condition to a
limit cycle in the macroscopic description is shown for $w=0.7$ in
\Fig{fig:trans}. 
The restriction operator $\restrict = (m,\dot m)$ is therefore defined
by the macroscopic description \eqref{eq:pedmacro} and its derivative.

The lifting operator $\lift$ uses information about the distribution
of the pedestrians in front of the door to initialize a sensible
microscopic state. The distribution of positions of pedestrians
  along the corridor is known from numerical studies and is observed
  to be well-approximated by a linear density distribution, \ie the distribution is
  of the form $p(|x|) = a|x| + b$, where $|x|$ is the distance from
  the door along the corridor axis. The slope $a$ and interception $b$
  are determined by simulations for all parameter values of
  interest. The lifting uses these distributions to map, \ie lift $(m,\dot m)$
  to a ``physically correct'' microscopic state. All velocities are
  initially set to $0$, such that we lift to a microscopic state with $\dot m = 0$
(see~\cite{Corradi2012} for details).

\subsection*{Numerical results}
\begin{figure}[t]
  \centering
  \begin{minipage}{0.52\textwidth}
 \includegraphics[width =\textwidth]{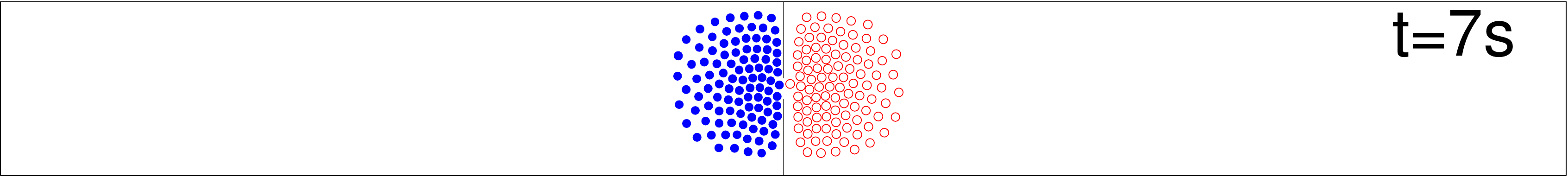}
\includegraphics[width =\textwidth]{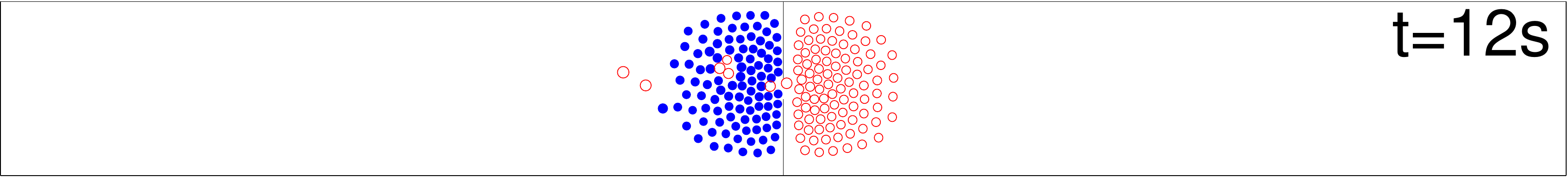} 
 \includegraphics[width =\textwidth]{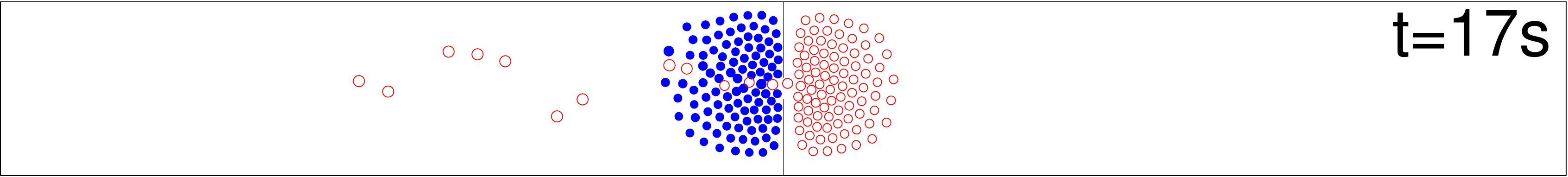}
\includegraphics[width =\textwidth]{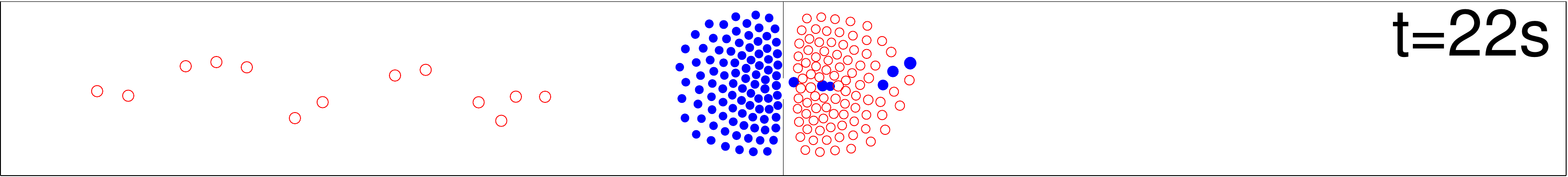} 
 \includegraphics[width =\textwidth]{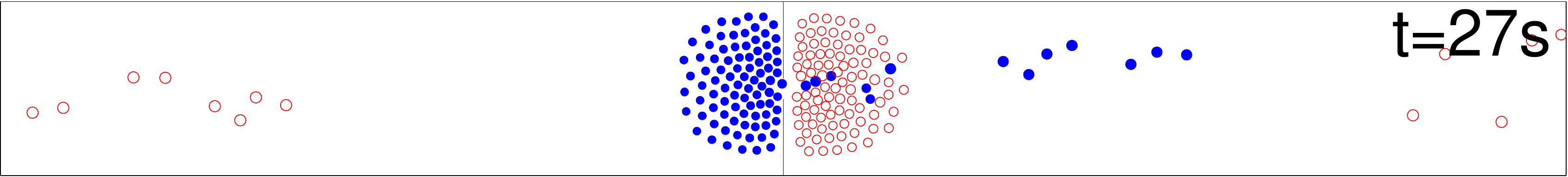}
\includegraphics[width =\textwidth]{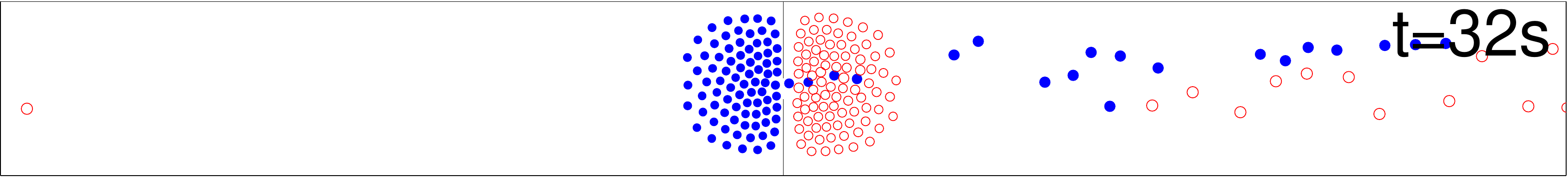} 
 \includegraphics[width =\textwidth]{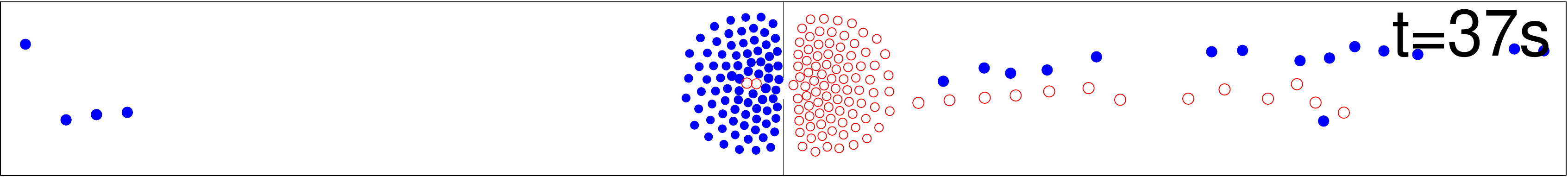}
\includegraphics[width =\textwidth]{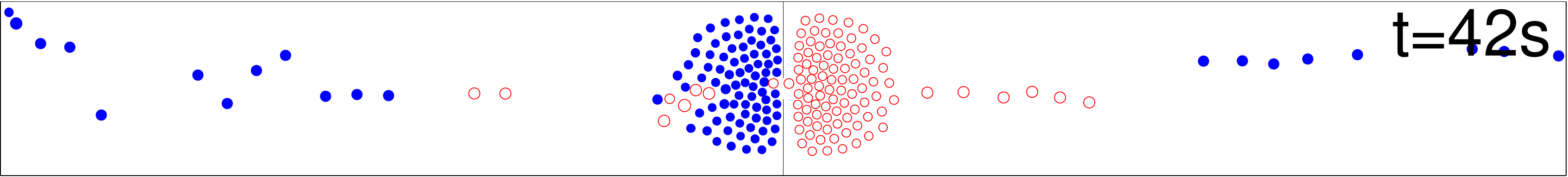} 
\subfigure[{Snapshots}]{\label{fig:snaps}\includegraphics[width
  =\textwidth]%
{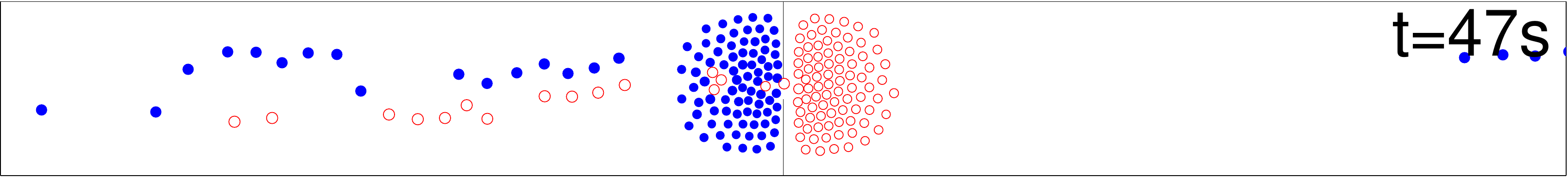}}
 \end{minipage}
  %
%
\begin{minipage}{0.38\textwidth}
\hfill
  \subfigure[{Two-parameter space}]{\label{fig:2ppeds}\includegraphics[width
    =.85\textwidth]{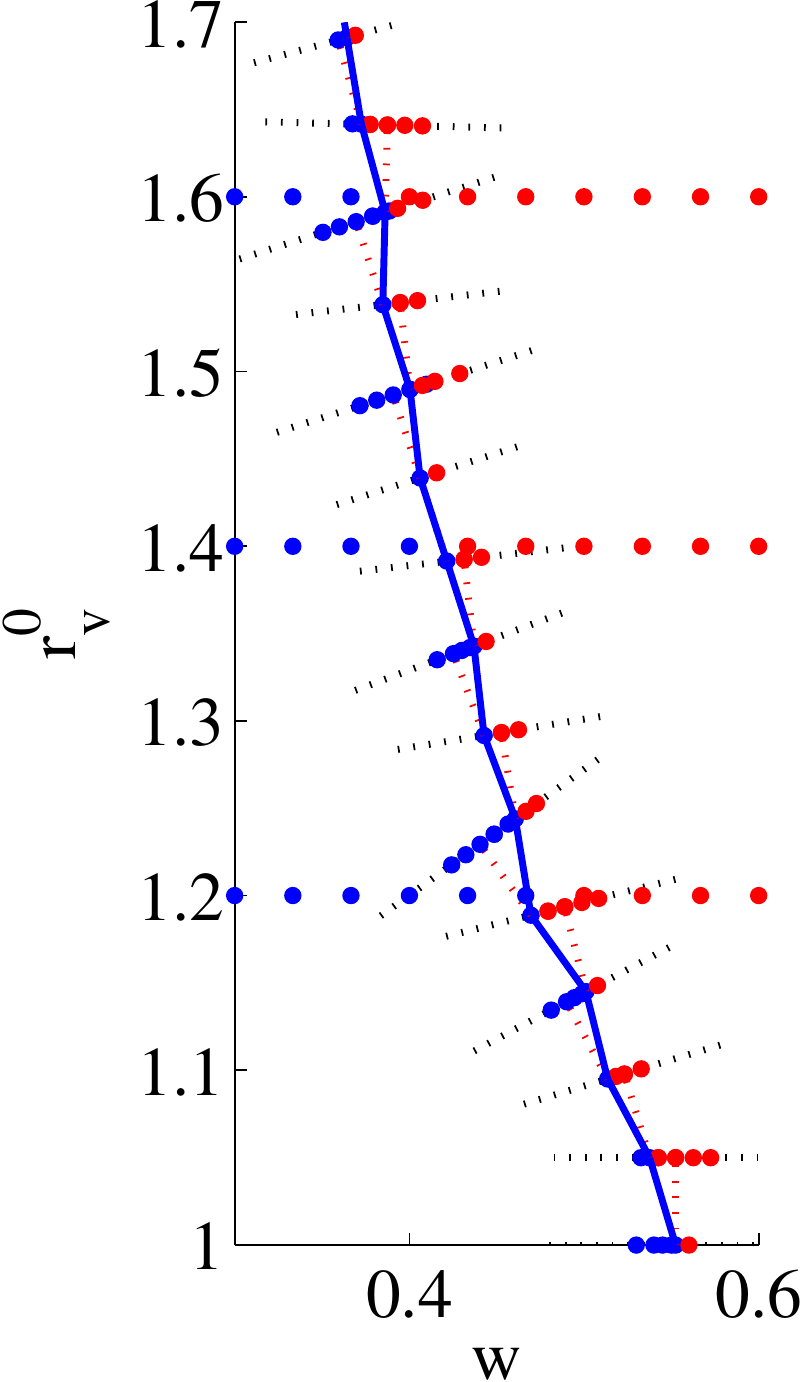}}%
\end{minipage}
\\
\subfigure[{Transient and limit cycle}]{\label{fig:trans}\includegraphics[width
  =.35\textwidth]%
{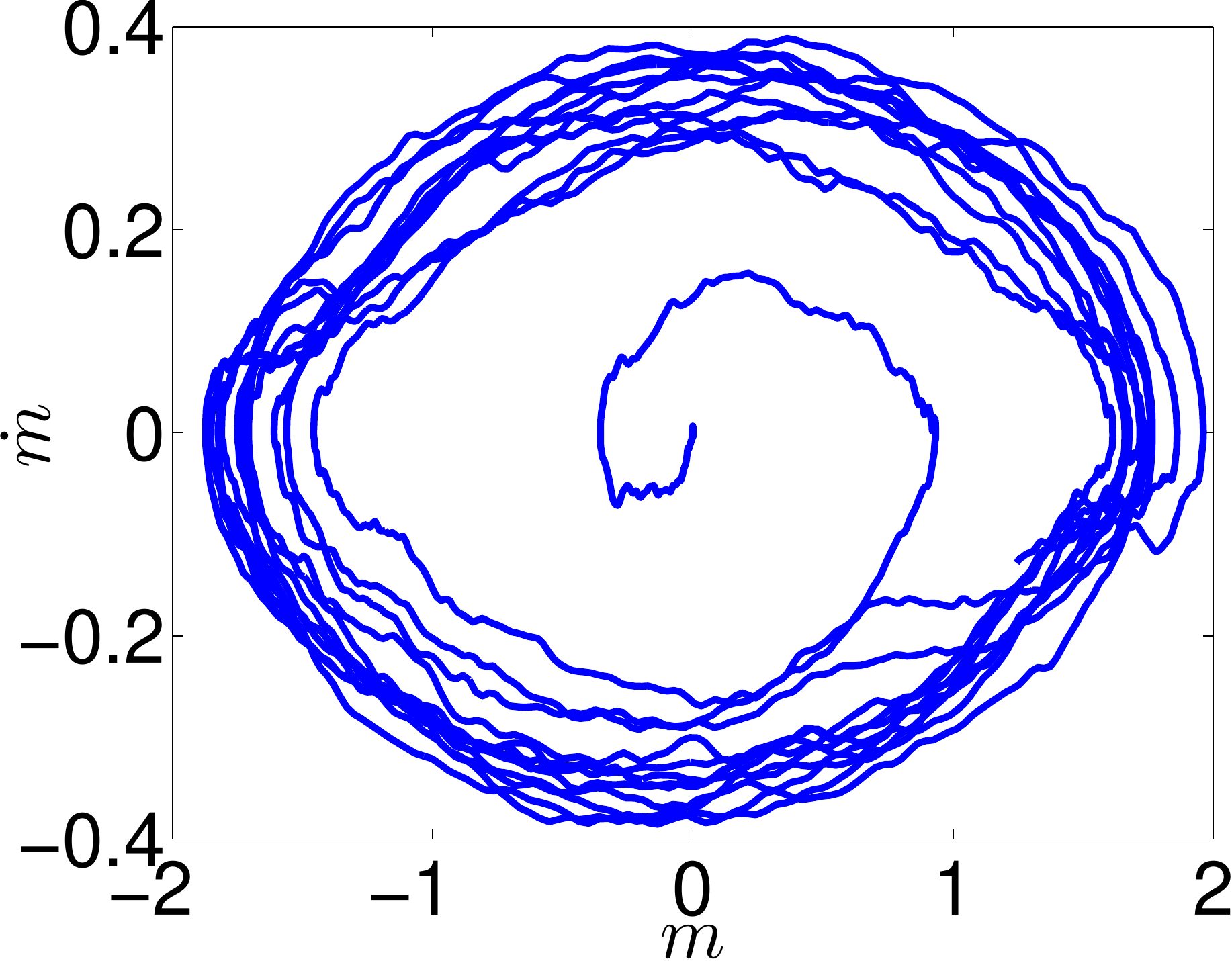}}
\hfill 
 \subfigure[{Bifurcation diagram}]{\label{fig:bifdiagped}\includegraphics[width =.35\textwidth]{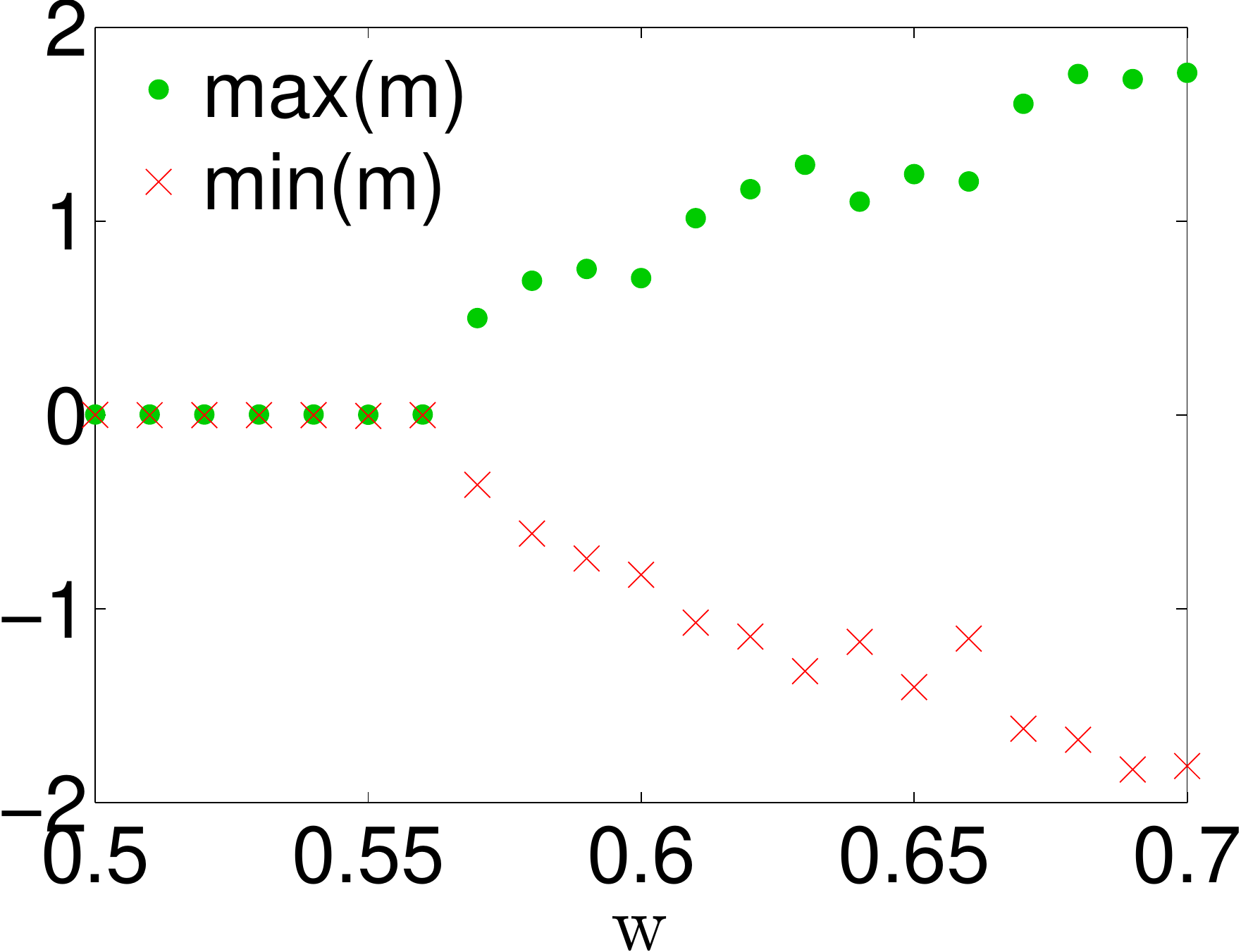}}
\caption{Coarse analysis of the pedestrian
  dynamics in a corridor with bottleneck. \textup{(a)} Snapshots of a
  microscopic simulation show oscillating behavior for large enough
  door width $w=0.6$. \textup{(b)}  Two-parameter plane explains the
  dynamics of the system and the point for the Hopf bifurcation. \textup{(c)} Transient and limit cycle in the
  macroscopic description for $w=0.7$. \textup{(d)} The coarse
  bifurcation diagram reveals a Hopf bifurcation at a critical door
  width $w=0.56$. }
  \label{fig:pedres}
\end{figure}
Using equation-free bifurcation analysis, the bifurcation diagram is
computed for the fixed ratio $r_{v_0} = 1$. \Fig{fig:bifdiagped} shows
the maximum and minimum of $m(t)$ as a function of $w$. The transition
from a blocked state to an oscillating state is clearly observed and
the bifurcation point is found to be at $w=0.56$. The transition is
analyzed in detail in~\cite{Corradi2012} and the bifurcation point is
identified as a Hopf bifurcation point using Poincar\'e sections, \ie
a discretization of the recurrent dynamics in time. This method is
also implicit with a healing time $\tskip$ determined
by the first crossing of the Poincar\'e section. The Hopf
bifurcation gives rise to macroscopic oscillations for large door
width $w$ emerging from a stable blocked state for $w$ small enough.

Let us now study the influence of $r_{v_0}$ on the location of the
bifurcation point. The system for macroscopic continuation is analyzed
by a predictor-corrector method using a linear prediction and a
subspace search for the correction in order to study the two-parameter
problem and to continue the Hopf point. The results are shown in
\Fig{fig:2ppeds}. Keeping the other model parameters fixed, this gives
an overview of the behavior of the system on a macroscopic level in
two parameters.

The application of equation-free analysis is not limited to
pedestrians in a bottleneck scenario. One could also think of
applications in evacuation scenarios (see, \eg
\cite{Helbing2002,Wagoum2012}), where parameter regimes with blocked
states have to be avoided at all cost. It is also possible to apply
equation-free analysis to discrete models, \eg cellular automaton
models \cite{Burstedde2001,Nowak2012}. This motivates further studies
using equation-free methods in traffic and pedestrian flow in order to
systematically investigate and finally optimize the parameter
dependencies of the macroscopic behavior of such microscopic models.

\section{Discussion and Conclusion}
\label{sec:conc}
We have demonstrated, that
equation-free methods can be useful to analyze the parameter dependent
behavior in traffic and pedestrian problems. Implicit methods allow us
to improve the results further by reducing the lifting error. The
comparison between traffic and pedestrian dynamics shows that both
problem classes can be studied with the same mathematical tools. In
particular, the use of coarse bifurcation analysis reveals some
information about the system that could not be obtained by
simpler means, \eg direct simulations of a microscopic model, since
they cannot investigate unstable solutions. Nevertheless, unstable
solutions are important in order to understand the phase space and parameter
dependence of the dynamics. In particular, in the case of a
one-dimensional macroscopic dynamics the unstable solutions act as
barriers between separate stable regimes defining reliable
  operating ranges. The knowledge of their
locations can be used to systematically push the system over the
barrier to switch to another more desirable solution, \eg leading
to a transition from traffic jams to uniform flow.  In the application
to two-dimensional macroscopic dynamics, we find the precise
  dividing line between oscillations and blocking in two
  parameters.

Finally, let us constrast equation-free analysis to the most obvious
alternative. A common approach to determining the precise parameter
value at which the onset of oscillations occurs, is to run the
simulation for sufficiently long time and observe if the transient
behavior vanishes.
This approach
suffers from two problems. First, close to the loss of linear
stability in the equilibrium (i.e.\ close to the bifurcation point)
the rate of approach to the stable orbit or fixed point is close to
zero as the Jacobi matrix becomes singular. This makes the transients
extremely long, resulting in unreliable numerics. Second, even
eventually decaying transients may grow intermittently (the effect of
\emph{non-normality}) such that the criteria for the choice of the
transient time to observe are non-trivial.  Equation-free computations
working on the macroscopic level in a neighborhood of the slow
manifold do not suffer from these long transients, as they are based
on direct root-finding methods.

In conclusion coarse bifurcation analysis can be used in future
research to improve safety in traffic problems and evacuation
scenarios of large buildings in case of emergency. The main advantage
is, that realistic models can be used and a qualitative analysis of
the macroscopic behavior is still possible. The method works almost
independent of the underlying microscopic model and has a significant
potential for helping traffic modellers to gain insight into
previously inaccessible scenarios.\\

\textbf{Acknowledgements} The authors thank their collaborators R.\
Berkemer, A.\ Kawamoto and O.\ Corradi. The research of
J.~Sieber is supported by EPSRC grant EP/J010820/1.
J.\ Starke was
partially funded by the Danish Research Council 
under 
09-065890/FTP and the Villum Fonden 
(VKR-Centre of
Excellence ``Ocean Life'').

%
%
\bibliographystyle{spphys}
\bibliography{tgfrefs}

\begin{thebibliography}{10}
\providecommand{\url}[1]{{#1}}
\providecommand{\urlprefix}{URL }
\expandafter\ifx\csname urlstyle\endcsname\relax
  \providecommand{\doi}[1]{DOI \discretionary{}{}{}#1}\else
  \providecommand{\doi}{DOI \discretionary{}{}{}\begingroup
  \urlstyle{rm}\Url}\fi

\bibitem{Helbing1995}
D.~Helbing, P.~Moln\'ar, Phys. Rev. E \textbf{51}, 4282 (1995)

\bibitem{Helbing2001}
D.~Helbing, Rev. Mod. Phys. \textbf{73}, 1067 (2001)

\bibitem{Corradi2012}
O.~Corradi, P.~Hjorth, J.~Starke, SIAM Journal on Applied Dynamical Systems
  \textbf{11}(3), 1007 (2012)

\bibitem{Bando1995}
M.~Bando, K.~Hasebe, A.~Nakayama, A.~Shibata, Y.~Sugiyama, Phys. Rev. E
  \textbf{51}(2), 1035 (1995)

\bibitem{Gasser2004}
I.~Gasser, G.~Sirito, B.~Werner, Physica D: Nonlinear Phenomena \textbf{197},
  222  (2004)

\bibitem{Sugiyama2008}
Y.~Sugiyama, M.~Fukui, M.~Kikuchi, K.~Hasebe, A.~Nakayama, K.~Nishinari, S.i.
  Tadaki, S.~Yukawa, New J. Phys. \textbf{10}(3) (2008)

\bibitem{Orosz2005}
G.~Orosz, B.~Krauskopf, R.~Wilson, Physica D: Nonlinear Phenomena \textbf{211},
  277  (2005)

\bibitem{Marschler2013}
C.~{Marschler}, J.~{Sieber}, R.~{Berkemer}, A.~{Kawamoto}, J.~{Starke}, ArXiv
  e-prints  (2013)

\bibitem{Haken1983a}
H.~Haken, \emph{Synergetics. An introduction.} (Springer, Berlin, 1983)

\bibitem{Haken1983}
H.~Haken, \emph{Advanced synergetics.} (Springer, Berlin, 1983)

\bibitem{Fenichel1979}
N.~Fenichel, Journal of Differential Equations \textbf{31}, 53 (1979)

\bibitem{kevrekidisgear2003}
I.G. Kevrekidis, C.W. Gear, J.M. Hyman, P.G. Kevrekidis, O.~Runborg,
  C.~Theodoropoulos, Communications in Mathematical Sciences \textbf{1}, 715
  (2003)

\bibitem{Kevrekidis2004}
I.G. Kevrekidis, C.W. Gear, G.~Hummer, AIChE Journal \textbf{50}(7), 1346
  (2004)

\bibitem{KevrekidisSamaey2009}
I.G. Kevrekidis, G.~Samaey, Annual Review of Physical Chemistry \textbf{60}(1),
  321 (2009)

\bibitem{Kevrekidis2010}
Y.~Kevrekidis, G.~Samaey, Scholarpedia \textbf{5}(9), 4847 (2010)

\bibitem{sieberkrauskopf08}
J.~Sieber, B.~Krauskopf, Nonlinear Dynamics \textbf{51}(3), 365 (2008)

\bibitem{Bureau2013}
E.~Bureau, F.~Schilder, I.F. Santos, J.J. Thomsen, J.~Starke, Journal of Sound
  and Vibration \textbf{332}(22), 5883  (2013)

\bibitem{Barton2013}
D.A.W. Barton, J.~Sieber, Phys. Rev. E \textbf{87}, 052916 (2013)

\bibitem{Gear2005}
C.W. Gear, T.J. Kaper, I.G. Kevrekidis, A.~Zagaris, SIAM Journal on Applied
  Dynamical Systems \textbf{4}, 711 (2005)

\bibitem{Zagaris2009}
A.~Zagaris, C.W. Gear, T.J. Kaper, Y.G. Kevrekidis, ESAIM: Mathematical
  Modelling and Numerical Analysis \textbf{43}(04), 757 (2009)

\bibitem{Zagaris2012}
A.~Zagaris, C.~Vandekerckhove, C.W. Gear, T.J. Kaper, I.G. Kevrekidis, Discrete
  and Continuous Dynamical Systems - Series A \textbf{32}(8), 2759  (2012)

\bibitem{Vandekerckhove2011}
C.~Vandekerckhove, B.~Sonday, A.~Makeev, D.~Roose, I.G. Kevrekidis, Computers
  \& Chemical Engineering \textbf{35}(10), 1949  (2011)

\bibitem{Beyn2002}
W.J. Beyn, A.~Champneys, E.~Doedel, W.~Govaerts, Y.A. Kuznetsov, B.~Sandstede,
  in \emph{Handbook of Dynamical Systems}, \emph{Handbook of Dynamical
  Systems}, vol.~2, ed. by B.~Fiedler (Elsevier Science, 2002), pp. 149 -- 219

\bibitem{Kuznetsov2004}
Y.A. Kuznetsov, \emph{Elements of Applied Bifurcation Theory}, \emph{Applied
  Mathematical Sciences}, vol. 112, 3rd edn. (Springer, New York, 2004)

\bibitem{Seyfried2006}
A.~Seyfried, B.~Steffen, T.~Lippert, Physica A: Statistical Mechanics and its
  Applications \textbf{368}(1), 232  (2006)

\bibitem{Zhang2012}
J.~Zhang, W.~Klingsch, A.~Schadschneider, A.~Seyfried, Journal of Statistical
  Mechanics: Theory and Experiment \textbf{2012}(02), P02002 (2012)

\bibitem{MarschlerStarkeLiuKevrekidis14}
C.~Marschler, J.~Starke, P.~Liu, I.G. Kevrekidis, Phys. Rev. E \textbf{89},
  013304 (2014)

\bibitem{Helbing2002}
D.~Helbing, I.J. Farkas, P.~Molnar, T.~Vicsek, Pedestrian and evacuation
  dynamics \textbf{21}, 21 (2002)

\bibitem{Wagoum2012}
A.U.K. Wagoum, M.~Chraibi, J.~Mehlich, A.~Seyfried, A.~Schadschneider, Computer
  Animation and Virtual Worlds \textbf{23}(1), 3 (2012)

\bibitem{Burstedde2001}
C.~Burstedde, K.~Klauck, A.~Schadschneider, J.~Zittartz, Physica A: Statistical
  Mechanics and its Applications \textbf{295}(3–4), 507  (2001)

\bibitem{Nowak2012}
S.~Nowak, A.~Schadschneider, Phys. Rev. E \textbf{85}, 066128 (2012)

\end{thebibliography}
\end{document}